\documentclass[12pt]{article}
\usepackage{amsfonts,amsmath,amssymb,euscript,amscd}
\usepackage{epic}
\usepackage {pstcol}
\usepackage{pstricks,pst-node,pst-coil,fancybox}
\usepackage[dvips]{graphics}
 \usepackage{graphicx}
 \newpsobject{showgrid}{psgrid}{subgriddiv=1,griddots=10,gridlabels=6pt}
\newcommand{\ov}{\over}
\newcommand{\iy}{\infty}

\newcommand{\s}{\sigma}
\newcommand{\g}{\gamma}

\newcommand{\bZ}{\mathbb{Z}}

\newcommand{\bP}{\mathbb{P}}
\renewcommand{\P}{\mathbb{P}}

\newcommand{\cI}{\mathcal{I}}

\newcommand{\be}{\begin{equation}}
\newcommand{\ee}{\end{equation}}
\newcommand{\bea}{\begin{eqnarray}}
\newcommand{\eea}{\end{eqnarray}}
\newcommand{\ba}{\begin{eqnarray*}}
\newcommand{\ea}{\end{eqnarray*}}
\newcommand{\bc}{\begin{center}}
\newcommand{\ec}{\end{center}}

\newcommand{\ra}{\rightarrow}

\renewcommand{\(}{\left(}
\renewcommand{\)}{\right)}

 \newcommand{\fr}{\frac}
\begin{document}
\hfill February 5, 2009
\begin{center}{\Large\bf Total Current Fluctuations in ASEP}\end{center}
 
\begin{center}{\large\bf Craig A.~Tracy}\\
{\it Department of Mathematics \\
University of California\\
Davis, CA 95616, USA\\
email: tracy@math.ucdavis.edu}\end{center}

\begin{center}{\large \bf Harold Widom}\\
{\it Department of Mathematics\\
University of California\\
Santa Cruz, CA 95064, USA\\
email: widom@ucsc.edu}\end{center}
\begin{center}\textbf{Abstract}\end{center}
A limit theorem for the total current in the asymmetric simple exclusion process (ASEP) with step initial condition is proved.  
This extends the result of Johansson on TASEP to ASEP. 
\section{Introduction}
The asymmetric simple exclusion process (ASEP) is a continuous time Markov process of interacting particles on the integer lattice $\bZ$ subject to two rules: (1) A particle at $x$ waits an exponential time with parameter one (independently of all other particles) and then it chooses $y$ with probability $p(x,y)$; (2) If $y$ is vacant at that time it moves to $y$, while if $y$ is occupied it remains at $x$ and restarts its clock.  The adjective ``simple'' refers to the fact that the allowed jumps are one step to the right, $p(x,x+1)=p$, or one step to the
left $p(x,x-1)=1-p=q$.  The asymmetric condition is $p\neq q$ so there is a net drift of particles.  
The special cases $p=1$ (particles hop only
to the right) or $q=1$ (particles hop only to the left) are called the T(totally)ASEP.
The dynamics are uniquely determined once we specify the initial 
state, which may be either deterministic
or random.  A rigorous construction of this infinite particle process can be found in Liggett \cite{liggett1}.

Since its introduction by Spitzer \cite{spitzer},  the ASEP has remained a popular model among probabilists  
 and physicists  because it is one of the simplest nontrivial processes modeling   nonequilibrium phenomena.  (For recent reviews see \cite{golinelli, liggett2,  seppalainen, spohn1}.)
If initially the particles are  located at $\bZ^{+}=\{1,2,\ldots\}$, called the \textit{step initial condition},
and if   $p<q$, then there will be on  average a net flow of particles, or \textit{current}, to the left.  More precisely,
we introduce the \textit{total current} $\cI$ at position $x\le0$ at time $t$:
\[\cI(x,t):=\# \>\> \textrm{of particles}\>\le x\>\>\textrm{at time}\>\> t.\]
With step initial condition, it has been known for some time (see, e.g., Theorem 5.12 in \cite{liggett1}) that  if we set
$\g:=q-p>0$ and $0\le c\le\g$, then the current $\cI$ satisfies the strong law
\[ \lim_{t\ra\iy}\fr{\cI([-c t],t)}{t}=\fr{1}{4\g}\(\g-c\)^2\ \ \textrm{a.s.}.\]
 
The natural next step is to examine the current fluctuations
\be \cI(x,t)-\fr{1}{4\g}(\g-c)^2\, t \label{diff}\ee
for large $x$ and $t$.
Physicists conjectured \cite{KPZ}, and Johansson proved for TASEP \cite{johansson},  that 
to obtain a nontrivial limiting distribution the correct  normalization of (\ref{diff})  is cube root in $t$.
For TASEP Johansson not only proved that the fluctuations are of order $t^{1/3}$ but found
 the limiting distribution
function.  Precisely, for $0\le v<1$  we have\footnote{The value of $a_2$ given in (\ref{a}) corrects a misprint in Corollary~1.7 of \cite{johansson}.}
\be \lim_{t\ra\iy}\bP\( \fr{\cI([-vt],t)-a_1 t }{ a_2 t^{1/3} }\le s\)=1-F_2(-s),\label{limitLaw}\ee
where
\be a_1=\fr{1}{4}\(1-v\)^2,\>\> a_2=2^{-4/3}(1-v^2)^{2/3},\label{a}\ee
and $F_2$ is the limiting distribution of the largest eigenvalue in the Gaussian Unitary Ensemble \cite{TW1}.

The proof of this relied  on the fact that TASEP is a determinantal process \cite{johansson,soshnikov,spohn1}.  However, universality arguments suggest that (\ref{limitLaw}) 
should extend to ASEP with step initial condition even though ASEP is not a determinantal process.  When the initial state is the Bernoulli product measure, it has been recently proved, using general probabilistic arguments, that the correct normalization remains $t^{1/3}$ for a large class of stochastic models including ASEP \cite{BKS, BS1, BS2, QV}.

In this paper we show that (\ref{limitLaw}) does extend to ASEP. 
\par\vspace{5mm}
\noindent\textbf{Theorem.} For ASEP with step initial condition we have, for $0\le v<1$,
\[\lim_{t\ra\iy}\bP\( \fr{\cI([-vt],t/\g)-a_1 t }{ a_2 t^{1/3} }\le s\)=1-F_2(-s),\]
where $\g=q-p$ and $a_1$ and $a_2$ are given by (\ref{a}).\footnote{With step initial condition and $x>0$ the total current equals the number of particles to the left of $x$ at time $t$ minus $x$. In what follows we shall require only that $|v|<1$. Therefore the statement of the Theorem holds for all such $v$ if when $v<0$ the value of $a_1$ is decreased by $|v|$.}
\par\vspace{5mm}
This theorem is a corollary, as we show below,  of 
earlier work by the authors \cite{TW2}.
\par\vspace{5mm}
\section{Proof of the Theorem}
We denote by $x_m(t)$ the position of the $m$th left-most particle (thus $x_m(0)=m\in\bZ^+$). 
We are interested in the probability of the event
\be \{\cI(x,t)=m\}=\{x_m(t)\le x, x_{m+1}(t)>x\}. \label{event}\ee
The sample space
consists of the four disjoint events $\{x_m(t)\le x, x_{m+1}(t)>x\}$,
 $\{x_m(t)\le x, x_{m+1}(t)\le x\}$,
  $\{x_m(t)> x, x_{m+1}(t)>x\}$,
   $\{x_m(t)\le x, x_{m+1}(t)\le x\}$ and because of the exclusion property we have
   \ba
   \{ x_m(t)\le x, x_{m+1}(t)\le x\}&=&\{x_{m+1}(t)\le x\},\\
   \{ x_m(t)>x,x_{m+1}(t)>x\}&=& \{x_m(t)>x\},\\
   \{x_m(t)>x,x_{m+1}(t)\le x\}&=& \emptyset.\ea
   These observations and (\ref{event}) give (the intuitively obvious)
   \[ \bP\(\cI(x,t)=m\)=\bP\(x_m(t)\le x\)-\bP\(x_{m+1}(t)\le x\).\]
   Since $\bP\(\cI(x,t)=0\)=\bP\(x_1(t)>x\)$, we have
   \[\bP\(\cI(x,t)\le m\)=1-\bP\(x_{m+1}(t)\le x\). \]
   Thus, since $x$ and $x_{m+1}(t)$ are integers, the statement of the Theorem is equivalent to the statement that
\[\lim_{t\to\iy}\P(x_{m+1}(t/\g)\le -vt)=F_2(s),\]
when $m=[a_1\,t-a_2\,s\,t^{1/3}]$.  
In fact, we shall show that
\be\lim_{t\to\iy}\P(x_{m}(t/\g)\le -vt)=F_2(s),\label{theq}\ee
when
\be m=a_1\,t-a_2\,s\,t^{1/3}+o(t^{1/3}).\label{m}\ee     
     
     Let 
   \[ \s=\fr{m}{t},\>\>\> c_1=-1+2\sqrt{\s},\>\>\> c_2=\s^{-1/6}(1-\sqrt{\s})^{2/3}.\]
   It was proved in \cite{TW2} that when $0\le p<q$,
 \be\lim_{t\to\iy}\P(x_m(t/\g)\le c_1\,t+s\,c_2\, t^{1/3})=F_2(s)\label{TWth}\ee
uniformly for $\s$ in a compact subset of $(0,\,1)$. 

To obtain (\ref{theq}) from this we determine $\s$ so that 
\[-vt=c_1\,t+s\,c_2\, t^{1/3}.\]
Thus,
\[v=1-2\sqrt\s-s\,\s^{-1/6}\,(1-\sqrt\s)^{2/3}\,t^{-2/3}.\]
Solving, we get 
\[\({1-v\ov2}\)^2=\s+s\,\s^{1/3}\,(1-\sqrt\s)^{2/3}\,t^{-2/3}+O\(t^{-4/3}\),\]
from which we deduce
\[\s=\({1-v\ov2}\)^2-s\,\({1-v\ov2}\)^{2/3}\,\({1+v\ov2}\)^{2/3}\,t^{-2/3}+O\(t^{-4/3}\)\]
\[=\({1-v\ov2}\)^2-s\,2^{-4/3}\,(1-v^2)^{2/3}\,t^{-2/3}+O\(t^{-4/3}\).\footnote{Since the condition on $\s$ is $0<\s<1$, the corresponding condition on $v$ is $|v|<1$, as was stated in footnote 2.}\]

By the uniformity of (\ref{TWth}) in $\s$ we get the same asymptotics if we replace the $\s$ we just computed by any $\s$ satisfying
\[\s=\({1-v\ov2}\)^2-s\,2^{-4/3}\,(1-v^2)^{2/3}\,t^{-2/3}+o(t^{-2/3}).\]
Since this is exactly the statement that $m=\s t$ satisfies (\ref{m}), we see that the Theorem is established.

\par\vspace{5mm}
\noindent\textbf{Acknowledgements.}  This work was supported by the National Science Foundation under grants
DMS--0553379 (first author) and DMS--0552388 (second author).


\begin{thebibliography}{99}

\bibitem{BKS} M.~Bal\'azs, J.~Komj\'athy and T.~Sepp\"al\"ainen, Microscopic concavity and
fluctuation bounds in a class of deposition processes, arXiv:0808.1177.

\bibitem{BS1} M.~Bal\'azs and T.~Sepp\"al\"ainen, Order of current variance and diffusivity in
the asymmetric simple exclusion process, arXiv:0806.0829, to appear in Annals of Math.

\bibitem{BS2}  M.~Bal\'azs and T.~Sepp\"al\"ainen, Fluctuation bounds for the asymmetric
simple exclusion process, arXiv:0806.0829.

\bibitem{golinelli} O.~Golinelli and K.~Mallick, The asymmetric simple exclusion process: an integrable
model for non-equilibrium statistical mechanics, J.\ Phys.\ A: Math.\ Gen.\ \textbf{39} (2006), 12679--12705.

\bibitem{johansson} K.~Johansson, Shape fluctuations and random matrices, Commun.\ Math.\ Phys.\
\textbf{209} (2000), 437--476.

\bibitem{KPZ} M.~Kardar, G.~Parisi and Y-C~Zhang, Dynamic scaling of growing interfaces,
Phys.\ Rev.\ Lett.\ \textbf{56} (1986), 889--892.

\bibitem{liggett1} T.~M.~Liggett, Interacting Particle Systems. Springer, Berlin (2005). (Reprint of the 1985 original.)
\bibitem{liggett2} T.~M.~Liggett, Stochastic Interacting Systems: Contact, Voter and Exclusion Processes.
Springer, Berlin (1999).

\bibitem{QV} J.~Quastel and B.~Valk\'o, $t^{1/3}$ superdiffusivity of finite-range asymmetric exclusion
processes on $\bZ$, Commun.\ Math.\ Phys.\ \textbf{273} (2007), 379--394.

\bibitem{seppalainen} T.~Sepp\"al\"ainen, Directed random growth models on the plane, arXiv:0708.2721.

\bibitem{soshnikov} A.~Soshnikov, Determinantal random fields, Russ.\ Math.\ Surv.\ \textbf{55} (2000), 923--975.
\bibitem{spitzer} F.~Spitzer, Interaction of Markov processes, Adv.\ Math.\ \textbf{5} (1970), 246--290.

\bibitem{spohn1} H.~Spohn, Exact solutions for KPZ-type growth processes, random matrices, and
equilibrium shapes of crystals, Phys.\ A \textbf{369} (2006), 71--99.

\bibitem{TW1} C.~A.~Tracy and H.~Widom, Level-spacing distributions and the Airy kernel,
Commun.\ Math.\ Phys.\ \textbf{159} (1994), 151--174.

\bibitem{TW2} C.~A.~Tracy and H.~Widom, Asymptotics in ASEP with step initial condition, arXiv:0807.1713.
\end{thebibliography}
\end{document}